\documentclass[a4paper,oneside,reqno,12pt]{amsart}
\usepackage{latexsym,amsthm,amscd,amsmath,amssymb}
\usepackage[mathscr]{eucal}

\def\bC{\mathbb C}
\def\bN{\mathbb N}
\def\bR{\mathbb R}

\def\bC{\mathbb C}
\def\cE{\mathcal E}
\def\cF{\mathcal F}
\def\cO{\mathcal O}
\def\cH{\mathcal H}

\def\cJ{\mathcal J}
\def\n{\noindent}

\def\PSH{\mathop{\rm PSH}\nolimits}
\def\interior{\mathop{\rm int}\nolimits}

\begin{document}
\centerline{\huge\bf A sharp lower bound for}
\vskip3mm
\centerline{\huge\bf the log canonical threshold}
\vskip 0.5cm
\n
\centerline{\bf Jean-Pierre Demailly and Ph\d{a}m Ho\`ang Hi\d{\^e}p}
\vskip3mm
\centerline{Institut Fourier, Universit\'e de Grenoble I,}
\centerline{Hanoi National University of Education}
\vskip 1cm

\n
{\bf Abstract.} In this note, we prove a sharp lower bound for the
log canonical threshold of a plurisubharmonic function $\varphi$ with
an isolated singularity at $0$ in an open subset of $\bC^n$. This
threshold is defined as the supremum of constants $c>0$ such that
$e^{-2c\varphi}$ is integrable on a neighborhood of~$0$. We relate
$c(\varphi)$ to the intermediate multiplicity numbers $e_j(\varphi)$, 
defined as the Lelong numbers of $(dd^c\varphi)^j$ at $0$ (so that
in particular $e_0(\varphi)=1$). Our main result is that
$c(\varphi)\ge\sum e_j(\varphi)/e_{j+1}(\varphi)$, $0\le j\le n-1$.
This inequality is shown to be sharp; it~simultaneously improves the
classical result $c(\varphi)\ge 1/e_1(\varphi)$ due to Skoda, as well
as the lower estimate $c(\varphi)\ge n/e_n(\varphi)^{1/n}$ which has
received crucial applications to birational geometry in recent
years. The proof consists in a reduction to the toric case, i.e.\
singularities arising from monomial ideals.

\vskip0.5cm
\n Keywords and Phrases:  Lelong number, Monge-Amp\`ere operator, log canonical threshold.

\n
\section{Notation and main results}

\n
Here we put $d^c=\frac {i} {2\pi} (\overline {\partial} - \partial)$, so that $dd^c=\frac {i} {\pi}\partial \overline {\partial}$. The normalization of the $d^c$ operator is chosen so that we have precisely $(dd^c\log|z|)^n=\delta_0$ for the Monge-Amp\`ere operator in~$\bC^n$. The Monge-Amp\`ere operator is defined on locally bounded plurisubharmonic functions according to the definition of Bedford-Taylor [BT76, BT82]; it can also be extended to plurisubharmonic functions with isolated or compactly supported poles by [Dem93]. 
If $\Omega$ is an open subset of $\bC^n$, we let $\PSH(\Omega)$ (resp.\ $\PSH^-(\Omega)$) be the set of plurisubharmonic (resp.\ psh${}\le 0$) functions on~$\Omega$.
\medskip

\n{\bf Definition 1.1.} {\it Let $\Omega$ be a bounded hyperconvex domain $($i.e.\ a domain possessing a negative psh exhaustion$)$. Following Cegrell {\rm [Ce04]}, 
we introduce certain classes of psh functions
on~$\Omega$, in relation with the definition of the Monge-Amp\`ere operator~$:$
$$\cE_0(\Omega)=\{\varphi\in \PSH^-(\Omega):\ \lim\limits_{z\to\partial\Omega}\varphi (z)=0,\ \int_\Omega (dd^c\varphi)^n<+\infty\},
\leqno{\rm(a)}$$
$$\cF(\Omega)=\{\varphi\in \PSH^-(\Omega):\ \exists\ \cE_0(\Omega)\ni\varphi_p\searrow\varphi,\ \sup\limits_{p\geq 1}\int_\Omega(dd^c\varphi_p)^n<+\infty\},\leqno{\rm(b)}$$
$$\cE(\Omega)=\{\varphi\in \PSH^-(\Omega):\ \exists\ \varphi_K\in\cF(\Omega)\ \text{such that}\ \varphi_K=\varphi\ \text{on}\ K,\ \forall K\subset\subset\Omega\}.\leqno{\rm(c)}$$}

\n It is proved in [Ce04] that the class $\cE(\Omega)$ is the biggest subset of $\PSH^-(\Omega)$ on which the Monge-Amp\`ere operator is well-defined. For a general
complex manifold $X$, after remo\-ving the negativity assumption of the functions involved,
one can in fact extend the Monge-Amp\`ere operator to the class 
$$\widetilde\cE(X)\subset\PSH(X)\leqno(1.2)$$ 
of psh functions which, on a neighborhood $\Omega\ni x_0$ of an arbitrary point 
$x_0\in X$, are
equal to a sum $u+v$ with $u\in\cE(\Omega)$ and $v\in C^\infty(\Omega)$; again, this is 
the biggest subclass of functions of $\PSH(X)$ on which the Monge-Amp\`ere operator
is locally well defined. It is easy to see that $\smash{\widetilde\cE(X)}$ contains the class of psh
functions which are locally bounded outside isolated singularities.\smallskip

For $\varphi\in\PSH(\Omega)$ and $0\in\Omega$, we introduce the log canonical threshold at $0$
$$c(\varphi) = \sup\big\{ c>0:\ e^{-2c\varphi} \text{ is } L^1 \text{ on a neighborhood of }0\big\},\leqno(1.3)$$
and for $\varphi\in\widetilde\cE(\Omega)$ we introduce the intersection numbers
$$e_{j} (\varphi ) = \int_{\{0\}} (dd^c\varphi)^j\wedge (dd^c\log\Vert z\Vert)^{n-j}
\leqno(1.4)$$
which can be seen also as the Lelong numbers of $(dd^c\varphi)^j$ at~$0$.
Our main result is the following sharp estimate. It is a generalization and a sharpening of similar inequalities discussed in [Cor95], [Cor00], [dFEM03], [dFEM04]; such inequalities have fundamental applications to birational geometry (see [IM72], [Puk87], [Puk02], [Isk01], [Che05]).
\medskip

\n
{\bf Theorem 1.5.} {\it Let $\varphi\in\widetilde\cE(\Omega)$ and $0\in\Omega$. Then $c(\varphi)=+\infty$
if $e_1(\varphi)=0$, and otherwise
$$c(\varphi)\geq\sum\limits_{j=0}^{n-1}\frac {e_j(\varphi)} {e_{j+1}(\varphi)}.$$}
\medskip

\n
{\bf Remark 1.6.} By Lemma 2.1 below, we have $(e_1(\varphi),\ldots,e_n(\varphi))\in D$
where
$$D=\big\{t=(t_1,\ldots,t_n)\in [0,+\infty)^n:\ t_1^2\leq t_2,\ t_{j}^2\leq t_{j-1}t_{j+1},\ \forall j=2,\ldots,n-1\big\},$$
i.e.\ $\log e_j(\varphi)$ is a convex sequence. In particular, we have 
$e_j(\varphi)\ge e_1(\varphi)^j$, and the denominators do not vanish in 1.5 if $e_1(\varphi)>0$. On the other hand, a well known inequality due to Skoda [Sko72] tells us that
$$\frac{1}{e_1(\varphi)}\leq c(\varphi)\leq \frac{n}{e_1(\varphi)},$$
hence $c(\varphi)<+\infty$ iff $e_1(\varphi)>0$.
To see that Theorem~1.5 is optimal, let us choose
$$\varphi (z) = \max\big(a_1\ln |z_1|,\ldots,a_n\ln |z_n|\big)$$
with $0<a_1\leq a_2\leq \ldots\leq a_n$. Then $e_j(\varphi)=a_1a_2\ldots a_j$, and a
change of variable $z_j=\zeta_j^{1/a_j}$ on $\bC\smallsetminus\bR_-$ easily shows that
$$c(\varphi)=\sum_{j=1}^n\frac{1}{a_j}.$$
Assume that we have a function $f:D\to [0,+\infty)$ such that 
$c(\varphi)\geq f (e_1(\varphi),\ldots,e_n(\varphi))$ for all 
$\varphi\in\widetilde\cE(\Omega)$. Then, by the above example, we must have
$$f(a_1,a_1a_2,\ldots,a_1\ldots a_n)\leq \sum\limits_{j=1}^n \frac 1 {a_j}$$
for all $a_j$ as above. By taking $a_j=t_j/t_{j-1}$, $t_0=1$, this implies that
$$f(t_1,\ldots,t_n)\leq \frac 1 {t_1}+\frac {t_1} {t_2}+\ldots+\frac {t_{n-1}} {t_n},\qquad\forall t\in D,$$
whence the optimality of our inequality. \qed
\medskip

\n
{\bf Remark 1.7.} Theorem~1.5 is of course stronger than Skoda's lower bound
$c(\varphi)\ge 1/e_1(\varphi)$. By the inequality between the arithmetic and
geometric means, we infer the main inequality of [dFEM03], [dFEM04] and [Dem09]
$$
c(\varphi)\ge \frac{n}{e_n(\varphi)^{1/n}}.\leqno(1.8)
$$
By applying the  arithmetic-geometric inequality for the indices $1\le j\le n-1$ in
our summation $\sum_{j=0}^{n-1}e_j(\varphi)/e_{j+1}(\varphi)$, we also infer the stronger inequality
$$c(\varphi)\geq \frac 1 { e_1(\varphi ) }+(n-1)\left[ \frac { e_1(\varphi ) } { e_n(\varphi ) }\right]^{\frac 1 {n-1}}.\leqno(1.9)$$
\vskip2mm

\section{Log convexity of the multiplicity sequence}

\n The log convexity of the multiplicity sequence can be derived from very
elementary integration by parts and the Cauchy-Schwarz inequality, using an argument
from [Ce04].\medskip

\n
{\bf Lemma 2.1}. {\it Let $\varphi\in\widetilde\cE(\Omega)$ and $0\in\Omega$. We have $e_{j}(\varphi )^2\leq e_{j-1} (\varphi ) e_{j+1} (\varphi ),\ \forall j=1,\ldots,n-1.$}
\medskip

\n
{\it Proof.} Without loss generality, by replacing $\varphi$ with a sequence of local 
approximations $\varphi_p(z)=\max(\varphi(z)-C,p\log|z|)$ of $\varphi(z)-C$, 
$C\gg 1$,  we can assume that $\Omega$ is the unit ball and $\varphi\in\cE_0(\Omega)$.
Take also $h,\,\psi\in\cE_0(\Omega)$. Then integration by parts and the Cauchy-Schwarz
inequality yield
$$\leqno\displaystyle
\left[\int_\Omega -h (dd^c\varphi)^j\wedge (dd^c\psi)^{n-j}\right]^2
=\left[\int_\Omega d\varphi\wedge d^c\psi \wedge (dd^c\varphi)^{j-1}\wedge (dd^c\psi)^{n-j-1}\wedge dd^ch\right]^2$$
$$\leq\int_\Omega d\psi\wedge d^c\psi \wedge (dd^c\varphi)^{j-1}\wedge (dd^c\psi)^{n-j-1}\wedge dd^ch\int_\Omega d \varphi\wedge d^c\varphi \wedge (dd^c\varphi)^{j-1}\wedge (dd^c\psi)^{n-j-1}\wedge dd^ch$$
$$=\int_\Omega -h (dd^c\varphi)^{j-1}\wedge (dd^c\psi)^{n-j+1}\int_\Omega -h (dd^c\varphi)^{j+1}\wedge (dd^c\psi)^{n-j-1}.$$
Now, as $p\to+\infty$, take
$$
h(z)=h_p(z)=\max\Big(-1,\frac{1}{p}\log\Vert z\Vert\Big)\nearrow \left\{
\begin{array}{ccc}
0 &\text{ if } z\in\Omega\smallsetminus \{0\}\\
-1 &\text{ if } z=0.\hfill
\end{array} \right.
$$
By the monotone convergence theorem we get in the limit
$$\left[\int_{ \{0\} } (dd^c\varphi)^j\wedge (dd^c\psi)^{n-j}\right]^2
\leq\int_{ \{0\} } (dd^c\varphi)^{j-1}\wedge (dd^c\psi)^{n-j+1}\int_{ \{0\} } 
(dd^c\varphi)^{j+1}\wedge (dd^c\psi)^{n-j-1}.$$
For $\psi(z)=\ln \Vert z\Vert$, this is the desired estimate.\qed
\medskip

\n
{\bf Corollary 2.2}. {\it Let $\varphi\in\widetilde\cE(\Omega)$ and $0\in\Omega$. We have the inequalities
\begin{eqnarray*}
e_j(\varphi)&\ge& e_1(\varphi)^j,\qquad \forall j=0,1,\ldots\le n\\
e_k(\varphi)&\le& e_j(\varphi)^{\frac{l-k}{l-j}}e_l(\varphi)^{\frac{k-j}{l-j}},
\qquad\forall 0\le j<k<l\le n.\\
\end{eqnarray*}
\vskip-3mm\n
In particular $e_1(\varphi)=0$ implies $e_k(\varphi)=0$ for $k=2,\ldots,n-1$ if $n\ge 3$.}
\medskip

\n
{\it Proof.} If $e_j(\varphi)>0$ for all $j$, Lemma 2.1 implies that $j\mapsto e_j(\varphi)/e_{j-1}(\varphi)$ is increasing, at least equal to $e_1(\varphi)/e_0(\varphi)=e_1(\varphi)$, and the inequalities follow from the log convexity. The general case can be proved by considering 
$\varphi_\varepsilon(z)=\varphi(z)+\varepsilon\log\Vert z\Vert$, since $0<\varepsilon^j\le 
e_j(\varphi_\varepsilon)\to
e_j(\varphi)$ when $\varepsilon\to 0$. The last statement is obtained by taking $j=1$ and~$l=n$.\qed
\vskip5mm

\section{Proof of the main theorem}
\n
We start with a monotonicity statement.\medskip

\n
{\bf Lemma 3.1.} {\it Let $\varphi,\psi\in\widetilde\cE(\Omega )$ be such that $\varphi \leq \psi$ $($i.e.\ $\varphi$ is ``more singular'' than $\psi)$. Then
$$\sum\limits_{j=0}^{n-1} \frac { e_{j} (\varphi ) } { e_{j+1} (\varphi ) } \leq \sum\limits_{j=0}^{n-1} \frac { e_{j} (\psi ) } { e_{j+1} (\psi ) }.$$}

\n{\it Proof.} As in Remark~1.6, we set
$$D=\{t=(t_1,\ldots,t_n)\in [0,+\infty)^n:\ t_1^2\leq t_2,\ t_{j}^2\leq t_{j-1}t_{j+1},\ \forall j=2,\ldots,n-1\}.$$
Then $D$ is a convex set in $\bR^n$, as can be checked by a straightforward application of
the Cauchy-Schwarz inequality. We consider the function $f: \interior D\to [0,+\infty)$ 
$$f(t_1,\ldots,t_n) = \frac 1 {t_1}+\frac {t_1} {t_2}\ldots+\frac {t_{n-1}} {t_n}.\leqno(3.2)$$
We have
$$\frac {\partial f} {\partial t_j} (t) = -\frac { t_{j-1} } { t_j^2 } + \frac { 1 } { t_{j+1} }\leq 0,\qquad \forall t\in D.$$
For $a,b\in\interior D$ such that $a_j\geq b_j$, $\forall j=1,\ldots,n$, $[0,1]\ni \lambda\to f( b+ \lambda(a-b) )$ is thus a decreasing function. This implies that $f(a)\leq f(b)$ for $a,b\in\interior D$, $a_j\geq b_j$, $\forall j=1,\ldots,n$. On the other hand, the hypothesis $\varphi\leq\psi$ implies $e_j (\varphi) \geq e_j (\psi)$, $\forall j=1,\ldots,n$, by the comparison principle (see e.g.\ [Dem87]). Therefore $f(e_1 (\varphi),\ldots,e_n (\varphi))\leq f(e_1 (\psi),\ldots,e_n (\psi))$.\qed
\medskip

\n
{\bf (3.3) Proof of the main theorem in the ``toric case''.}
\smallskip

\n
It will be convenient here to introduce Kiselman's refined Lelong numbers 
(cf.\ [Kis87], [Kis94a]):
\medskip

\n
{\bf Definition 3.4.} {\it Let $\varphi\in\PSH(\Omega)$. Then the function
$$\nu_{\varphi}(x)=\lim\limits_{t\to -\infty} \frac {\max\{\varphi (z): |z_1|=e^{x_1t},\ldots,|z_n|=e^{x_nt}\}} {t}$$
is called the refined Lelong number of $\varphi$ at $0$. This function is increasing in each variable $x_j$ and concave on $\bR_+^n$.}
\medskip

\n
By ``toric case'', we mean that $\varphi (z_1,\ldots,z_n) =
\varphi (|z_1|,\ldots,|z_n|)$ depends only on $|z_j|$ for all~$j$; then 
$\varphi$ is psh if and only if $(t_1,\ldots,t_n)\mapsto\varphi(e^{t_1},
\ldots,e^{t_n})$ is increasing in each $t_j$ and convex.
By replacing $\varphi$ with $\varphi(\lambda z)-\varphi(\lambda,...,\lambda)$, $0<\lambda\ll 1$, we can assume that $\Omega = \Delta^n$ is the unit polydisk, $\varphi(1,\ldots,1)=0$ (so that 
$\varphi\le 0$ on $\Omega$), and we have $e_1 (\varphi) = n\,\nu_\varphi(\frac 1 n,\ldots,\frac 1 n)$. By convexity, the slope $\frac{\max\{\varphi(z)\,:\;|z_j|=e^{x_jt}\}}{t}$ is 
increasing in $t$ for $t<0$. Therefore, by taking $t=-1$ we get
$$
\nu_\varphi(-\ln |z_1|,\ldots,-\ln |z_n|)\le -\varphi(z_1,\ldots,z_n).
$$
Notice also that $\nu_\varphi(x)$ satisfies the $1$-homogeneity property
$\nu_\varphi(\lambda x)=\lambda\nu_\varphi(x)$ for $\lambda\in\bR_+$. As~a
consequence, $\nu_\varphi$ is entirely characterized by its restriction to the set
$$\Sigma=\Big\{ x=(x_1,\ldots,x_n)\in\bR_+^n:\ \sum\limits_{j=1}^n x_j = 1 \Big\}.$$
We choose $x^0=(x_1^0,\ldots,x_n^0)\in\Sigma$ such that
$$\nu_{\varphi} (x^0) = \max\{\nu_{\varphi} (x):\ x\in \Sigma\}\in \Big[\frac{e_1(\varphi)}{n},e_1(\varphi)\Big].$$
By Theorem 5.8 in [Kis94a] (see also [Ho01] for similar results in an algebraic\
 context) we have the formula
$$c(\varphi) = \frac 1 { \nu_\varphi (x^0) }.$$
Set
$$\zeta (x) = \nu_\varphi(x^0) \min\Big(\frac { x_1 } {x_1^0},\ldots,\frac { x_n } {x_n^0}\Big),\qquad \forall x\in\bR_+^n.$$
Then $\zeta$ is the smallest nonnegative concave $1$-homogeneous function on $\bR_+^n$ that is increasing in each variable $x_j$ and such that $\zeta (x^0)=\nu_\varphi (x^0)$. Therefore we have $\zeta\le \nu_\varphi$, hence
\begin{eqnarray*}
\varphi (z_1,\ldots,z_n) &\leq& -\nu_{\varphi} (-\ln |z_1|,\ldots,-\ln |z_n|)\\
&\leq& -\zeta (-\ln |z_1|,\ldots,-\ln |z_n|)\\
&\leq& \nu_{\varphi} (x^0)\max\left( \frac { \ln |z_1| } {x_1^0},\ldots,\frac { \ln |z_n| } {x_n^0} \right):=\psi (z_1,\ldots,z_n).
\end{eqnarray*}
By Lemma 3.1 and Remark 1.6 we get
$$f(e_1 (\varphi),\ldots,e_n (\varphi))\leq f(e_1 (\psi),\ldots,e_n (\psi))=c(\psi)=\frac{1}{\nu_\varphi(x^0)}=c(\varphi).$$
\vskip3pt

\n
{\bf (3.5) Reduction to the case of psh functions with analytic singularities.}

\n
In the second step, we reduce the proof to the case $\varphi = \log (|f_1|^2+\ldots+|f_N|^2)$, where $f_1,\ldots,f_N$ are germs of holomorphic functions at $0$. Following the technique introduced in~[Dem92], we let $\cH_{m\varphi}(\Omega)$ be the Hilbert space of holomorphic functions $f$ on $\Omega$ such that
$$\int_\Omega|f|^2e^{-2m\varphi}dV<+\infty,$$
and let $\psi_m=\frac{1}{2m}\log\sum|g_{m,k}|^2$ where $\{g_{m,k}\}_{k\geq 1}$ is an orthonormal basis of $\cH_{m\varphi}(\Omega)$. Thanks to Theorem 4.2 in [DK00], mainly based on to the Ohsawa-Takegoshi $L^2$ extension theorem [OT87] (see also [Dem92]), there are constants $C_1,C_2>0$ independent of $m$ such that
$$
\varphi(z)-\frac{C_1}{m}\le
\psi_m(z)\le\sup_{|\zeta-z|<r}\varphi(\zeta)+\frac{1}{m}\log\frac{C_2}{r^n}
$$
for every $z\in\Omega$ and $r<d(z,\partial\Omega)$ and
$$\nu(\varphi)-\frac{n}{m}\le\nu(\psi_m)\le\nu(\varphi),$$
$$\frac 1 {c(\varphi)}-\frac{1}{m}\le\frac{1}{c(\psi_m)}\le \frac{1}{c(\varphi)}.$$
By Lemma 3.1, we have
$$f(e_1 (\varphi),\ldots,e_n (\varphi))\leq f(e_1 (\psi_m),\ldots,e_n (\psi_m)),\qquad \forall m\geq 1.$$
The above inequalities show that in order to prove the lower bound of $c(\varphi)$ in Theorem~1.5, we only need to prove it for $c(\psi_m)$ and let $m$ tend to infinity. Also notice that since the Lelong numbers of a function $\varphi\in\smash{\widetilde\cE}(\Omega)$ occur only on a discrete set, the same is true for the functions~$\psi_m$.
\medskip

\n
{\bf (3.6) Reduction of the main theorem to the case of monomial ideals.}
\smallskip

\n
The final step consists of proving the theorem for $\varphi = \log (|f_1|^2+\ldots.+|f_N|^2)$, where $f_1,\ldots,f_N$ are germs of holomorphic functions at $0$ [this is because the ideals $(g_{m,k})_{k\in\bN}$ in the Noetherian ring $\cO_{\bC^n,0}$ are always finitely generated]. Set $\cJ=(f_1,\ldots,f_N)$, $c(\cJ) = c(\varphi)$, $e_j(\cJ) = e_j(\varphi)$, $\forall j=0,\ldots,n$. By the final observation of 3.5, we can assume that $\cJ$ has an isolated zero at~$0$. Now, by fixing a multiplicative order on the monomials $z^\alpha=z_1^{\alpha_1}\ldots z_n^{\alpha_n}$ (see [Eis95] Chap.~15 and [dFEM04]), it is well known that one can construct a flat family $(\cJ_s)_{s\in\bC}$ of ideals of $\cO_{\bC^n,0}$ depending on a complex parameter $s\in\bC$, such that $\cJ_0$ is a monomial ideal, $\cJ_1=\cJ$ and $\dim(\cO_{\bC^n,0}/\cJ_s^t)=\dim(\cO_{\bC^n,0}/\cJ^t)$ for all $s$ and $t\in\bN$; in fact $\cJ_0$ is just the initial ideal associated to $\cJ$ with respect to the monomial order. Moreover, we can arrange by a generic rotation of
coordinates $\bC^p\subset\bC^n$ that the family of ideals 
$\cJ_{s\,|\,\bC^p}$ is also flat, and that the dimensions
$$
\dim(\cO_{\bC^p,0}/(\cJ_{s\,|\,\bC^p})^t)=\dim(\cO_{\bC^p,0}/(\cJ_{\,|\,\bC^p})^t)
$$
compute the intermediate multiplicities
$$
e_p(\cJ_s)=\lim_{t\to+\infty}\frac{p!}{t^p}\dim(\cO_{\bC^p,0}/(\cJ_{s\,|\,\bC^p})^t)=e_p(\cJ)
$$
(notice, in the analytic setting, that the Lelong number of the
$(p,p)$-current $(dd^c\varphi)^p$ at $0$ is the Lelong number of
its slice on a generic $\bC^p\subset\bC^n$); in particular $e_p(\cJ_0)=e_p(\cJ)$ for all~$p$. The semicontinuity property of the log canonical threshold (see for example [DK00]) now implies that $c(\cJ_0)\leq c(\cJ_s)$ for $s$ small.
As $c(\cJ_s)=c(\cJ)$ for $s\ne 0$ ($\cJ_s$ being a pull-back of $\cJ$ by a biholomorphism, in other words $\cO_{\bC^n,0}/\cJ_s\simeq \cO_{\bC^n,0}/\cJ$ as rings,
see again [Eis95], chap.~15), the lower bound is valid for $c(\cJ)$ if it is valid for $c(\cJ_0)$.

\vskip 0.5cm
\n
Jean-Pierre Demailly

\n
Universit\'e de Grenoble I, D\'epartement de Math\'ematiques

\n
Institut Fourier, 38402 Saint-Martin d'H\`eres, France

\n
{\it e-mail\/}: {\tt jean-pierre.demailly@ujf-grenoble.fr}
\vskip 0.5 cm
\n
Ph\d{a}m Ho\`ang Hi\d{\^e}p

\n
Department of Mathematics, National University of Education

\n
136-Xuan Thuy, Cau Giay, Hanoi, Vietnam

\n 
and Institut Fourier (on a Post-Doctoral grant from Univ.\ Grenoble I)

\n
{\it e-mail\/}: {\tt phhiep$_-$vn@yahoo.com}

\vskip1cm

\begin{thebibliography}{00000000}

%% \bibitem[ACCH]{} P.\ \AA hag, U.\ Cegrell, R.\ Czyz, H.\ H.\ Ph\d{a}m, {\it Monge-Amp\`{e}re measure on pluripolar sets}, J.\ Math.\ Pures Appl.\ {\bf 92} (2009), 613--627.

%% \bibitem[ACKHZ]{} P.\ \AA hag, U.\ Cegrell, S.\ Ko{\l}odziej, Ph\d{a}m Ho\`ang Hi\d{\^e}p, A.\ Zeriahi, {\it Partial pluricomplex energy and integrability exponents of plurisubharmonic functions}, Adv.\ Math.\ {\bf 222} (2009), 2036--2058.

%% \bibitem[BAC83]{} N.\ Bourbaki, {\it Alg\`ebre Commutative, chapter VIII et IX}; Masson, Paris, 1983.

\bibitem[BT76]{} E.\ Bedford and B.\ A.\ Taylor, {\it The Dirichlet problem for a complex Monge-Amp\`ere equation}, Invent.\ Math.\ {\bf 37} (1976) 1--44.

\bibitem[BT82]{} E.\ Bedford and B.\ A.\ Taylor, {\it A new capacity for plurisubharmonic functions}, Acta Math.\ {\bf 149} (1982) 1--41.

\bibitem [Ce04]{} U.\ Cegrell, {\it The general definition of the complex Monge-Amp\`ere operator}, Ann.\ Inst.\ Fourier (Grenoble)., {\bf 54}(2004), 159--179.

\bibitem[Che05]{} I.\ Chel'tsov, {\it Birationally rigid Fano manifolds}, Uspekhi Mat.\ Nauk {\bf 60:}5 (2005), 71-160 and Russian Math.\ Surveys {\bf 60:}5 (2005), 875--965.

\bibitem[Cor95]{} A.\ Corti, {\it Factoring birational maps of threefolds after Sarkisov}, J.\ Algebraic Geom., {\bf 4} (1995), 223--254.

\bibitem[Cor00]{} A.\ Corti, {\it Singularities of linear systems and $3$-fold birational geometry}, London Math.\ Soc.\ Lecture Note Ser.\ { \bf 281} (2000) 259--312.

\bibitem[dFEM03]{} T.\ de Fernex, T, L.\ Ein and Musta\c{t}\v{a}, {\it Bounds for log canonical thresholds with applications to birational rigidity}, Math.\ Res.\ Lett.\ {\bf 10} (2003) 219--236.

\bibitem[dFEM04]{} T.\ de Fernex, T, L.\ Ein and Musta\c{t}\v{a}, {\it Multiplicities and log canonical thresholds}, J.\ Algebraic Geom.\ {\bf 13} (2004) 603--615.

\bibitem[Dem87]{} J.-P.\ Demailly, {\it Nombres de Lelong g\'en\'eralis\'es, th\'eor\`emes d'int\'egralit\'e et d'analyticit\'e}, Acta Math.\ {\bf 159} (1987) 153--169.

%% \bibitem[Dem90]{} J.-P.\ Demailly, {\it Singular hermitian metrics on positive line bundles}, Proceedings of the Bayreuth conference ``Complex algebraic varieties'', April 2-6, 1990, edited by K.\ Hulek, T.\ Peternell, M.\ Schneider, F.\ Schreyer, Lecture Notes in Math.\ n${}^\circ\,$1507, Springer-Verlag, 1992.

\bibitem[Dem92]{} J.-P.\ Demailly, {\it Regularization of closed positive currents and Intersection Theory}, J.\ Alg.\ Geom.\ {\bf 1} (1992), 361--409.

\bibitem[Dem93]{} J.-P.\ Demailly, {\it Monge-Amp\`ere operators, Lelong numbers and intersection theory}, Complex Analysis and Geometry, Univ.\ Series in Math., edited by V.\ Ancona and A.\ Silva, Plenum Press, New-York, 1993.

\bibitem[DK00]{} J.-P.\ Demailly and J.\ Koll\'ar, {\it Semi-continuity of complex singularity exponents and K\"ahler-Einstein metrics on Fano orbifolds}, Ann.\ Sci.\ Ecole Norm.\ Sup.\ (4) {\bf 34} (2001), 525--556.

\bibitem[Dem09]{} J.-P.\ Demailly, {\it Estimates on Monge-Amp\`ere operators derived from a local algebra inequality}, in: Complex Analysis and Digital geometry, Proceedings of the Kiselmanfest 2006, Acta Universitatis Upsaliensis, 2009.

\bibitem[Eis95]{} D.\ Eisenbud, {\it Commutative algebra with a view toward algebraic geometry}, Grad.\ Texts in Math.\ {\bf 150}, Springer, New York, 1995.

\bibitem[Ho01]{} J.\ Howald, {\it Multiplier ideals of monomial ideals}, Trans.\ Amer.\ Math.\ Soc.\ {\bf 353} (2001), 2665--2671.

\bibitem[IM72]{} V.A.\ Iskovskikh and I.Yu.\ Manin, {\it Three-dimensional quartics and counterexamples to the L\"uroth problem}, Mat.\ Sb.\ {\bf 86} (1971), 140--166; English transl., Math.\ Sb.\ {\bf 15} (1972), 141--166.

\bibitem[Isk01]{} V.A.\ Iskovskikh, {\it Birational rigidity and Mori theory}, Uspekhi Mat.\ Nauk {\bf 56}:2 (2001) 3-86, English transl., Russian Math.\ Surveys {\bf 56}:2 (2001), 207--291.

%% \bibitem[Kis78]{} C.O.\ Kiselman, {\it The partial Legendre transformation for plurisubharmonic functions}, Invent.\ Math.\ {\bf 49} (1978) 137--148.

%% \bibitem[Kis79]{} C.O.\ Kiselman, {\it Densit\'e des fonctions plurisousharmoniques}, Bulletin de la Soci\'et\'e Math\'ematique de France, {\bf 107} (1979) 295--304.

%% \bibitem[Kis84]{} C.O.\ Kiselman, {\it Sur la d\'efinition de l'op\'erateur de Monge-Amp\`ere complexe}, Analyse Complexe, Proceedings of the Journ\'ees Fermat -- Journ\'ees SMF, Toulouse 1983; Lecture Notes in Mathematics {\bf 1094}, Springer-Verlag (1984) 139--150.

\bibitem[Kis87]{} C.O.\ Kiselman, {\it Un nombre de Lelong raffin\'e}, S\'eminaire d'Analyse Complexe et G\'eom\'etrie 1985­87, Facult\'e des Sciences de Tunis \& Facult\'e des Sciences et Techniques de Monastir, (1987) 61--70.

\bibitem[Kis94a]{} C.O.\ Kiselman, {\it Attenuating the singularities of plurisubharmonic functions}, Ann.\ Polon.\ Math.\ {\bf 60} (1994) 173--197.

%% \bibitem[Kis94b]{} C.O.\ Kiselman, {\it Plurisubharmonic functions and their singularities}, Complex Potential Theory (Eds.\ P.\ M.\ Gauthier G.\ Sabidussi).\ NATO ASI Series, Series C Vol.\ {\bf 439} Kluwer Academic Publishers (1994) 273--323.

\bibitem[OhT87]{} T.\ Ohsawa and K.\ Takegoshi, {\it On the extension of $L^2$ holomorphic functions}, Math.\ Zeit.\ {\bf 195} (1987)
  197--204.

\bibitem[Puk87]{} A.V.\ Pukhlikov, {\it Birational automorphisms of a four-dimensional quintic}, Invent.\ Math.\ {\bf 87} (1987), 303--329.

\bibitem[Puk02]{} A.V.\ Pukhlikov, {\it Birationally rigid Fano hypersurfaces}, Izv.\ Ross.\ Akad.\ Nauk Ser.\ Mat.\ {\bf 66:}6 (2002), 159-186; English translation, Izv.\ Math.\ {\bf 66:} (2002), 1243--1269.

\bibitem[Sko72]{} H.\ Skoda, {\it Sous-ensembles analytiques d'ordre fini ou infini dans $\bC^n$}, Bull.\ Soc.\ Math.\ France {\bf 100} (1972) 353--408.

\end{thebibliography}
\end{document}